\documentclass{article}
\usepackage{latexsym}\usepackage{amsbsy}\usepackage{amssymb}
\usepackage[latin1]{inputenc}

\newtheorem{theorem}{Theorem}[section]
\newtheorem{lemma}[theorem]{Lemma}
\newtheorem{proposition}[theorem]{Proposition}

\newtheorem{definition}[theorem]{Definition}
\newtheorem{example}[theorem]{Example}

\newtheorem{remark}[theorem]{Remark}

\newcommand{\pro}{\mathop{\mathrm{pr}}}
\newcommand{\eps}{\varepsilon}
\newcommand{\ms}{\medskip\\}
\newcommand{\cl}{}

\newcommand{\R}{\mathbb R}
\newcommand{\N}{\mathbb N}
\newcommand{\C}{\mathbb C}

\newcommand{\K}{\mathbb K}

\newcommand{\gK}{{\mathcal K}} 
\newcommand{\gR}{{\mathcal R}} 

\newcommand{\gs}{\ensuremath{{\mathcal G}} }

\newcommand{\esm}{\ensuremath{{\mathcal E}_M} }
\newcommand{\esmh}{\ensuremath{{\mathcal E}^{\mathit h}_M} }
\newcommand{\gsh}{\ensuremath{{\mathcal G}^{\mathit h}} }
\newcommand{\simh}{\ensuremath{\sim_{\mathit h}}}

\newcommand{\ns}{\ensuremath{{\mathcal N}} }

\newcommand{\ks}{\ensuremath{{\mathcal K}} }

\newcommand{\comp}{\subset\subset}

\newcommand{\cinfty}{{\mathcal C}^\infty}

\newcommand{\beast}{\begin{eqnarray*}}
\newcommand{\eeast}{\end{eqnarray*}}

\newcommand{\rem}[1]{\vadjust{\rlap{\kern\hsize\thinspace\vbox%
                       to0pt{\hbox{${}_\clubsuit${\small\tt #1}}\vss}}}}

\newcommand{\al}{\alpha}
\newcommand{\bet}{\beta} 
\newcommand{\ga}{\gamma}
\newcommand{\Om}{\Omega}\newcommand{\Ga}{\Gamma}
\newcommand{\la}{\lambda}
\newcommand{\de}{\delta}
\newcommand{\vphi}{\varphi}

\newcommand{\D}{{\mathcal D}}\newcommand{\Vol}{\mbox{Vol\,}}

\newcommand{\pa}{\partial}
\newcommand{\T}{{\mathcal T}} \newcommand{\G}{{\mathcal G}}
 
\newcommand{\CC}{{\mathcal C}}

\newcommand{\nn}{\nonumber}

\newcommand{\be}{ \begin{equation} }\newcommand{\ee}{\end{equation} }
\newcommand{\beq}{ \begin{equation} }\newcommand{\eeq}{\end{equation} }
\newcommand{\bea}{\begin{eqnarray}}\newcommand{\eea}{\end{eqnarray}}
\newcommand{\beas}{\begin{eqnarray*}}\newcommand{\eeas}{\end{eqnarray*}}
\newcommand{\beqs}{\begin{equation*}}\newcommand{\eeqs}{\end{equation*}}

\newcommand{\bfs}{\boldsymbol}
\newcommand{\ben}{\begin{enumerate}}\newcommand{\een}{\end{enumerate}}
\newcommand{\ba}{\begin{array}}\newcommand{\ea}{\end{array}}

\newcommand{\cc}{{\mathcal C}}

\newcommand{\ep}{\hspace*{\fill}$\Box$}
\newcommand{\pr}{{\bf Proof. }}
\newcommand{\bt}{\begin{theorem}}
\newcommand{\bp}{\begin{proposition}}
\newcommand{\bc}{\begin{corolloray}}
\newcommand{\blem}{\begin{lemma}}
\newcommand{\bex}{\begin{example}}
\newcommand{\bd}{\begin{definition}}
\newcommand{\brem}{\begin{remark}}
\begin{document}
\title{Generalized pseudo-Riemannian geometry}
\author{Michael Kunzinger \footnote{Electronic mail: michael.kunzinger@univie.ac.at}\\ 
        Roland Steinbauer \footnote{Electronic mail: roland.steinbauer@univie.ac.at}\\ 
         {\small Department of Mathematics, University of Vienna}\\
         {\small Strudlhofg.\ 4, A-1090 Wien, Austria}        }
\date{}  
\maketitle
\begin{abstract}
Generalized tensor analysis in the sense of Colombeau's
construction is employed to introduce a nonlinear distributional pseudo-Riemannian 
geometry. In particular, after deriving several characterizations of invertibility in 
the algebra of generalized functions we define the notions of generalized 
pseudo-Riemannian metric, generalized connection and generalized curvature 
tensor. We prove a ``Fundamental Lemma of \hbox{(pseudo-)} Riemannian geometry'' 
in this setting and define the notion of geodesics of a generalized metric. Finally, we  
present applications of the resulting theory to general relativity.
 
\end{abstract}

\section{Introduction}\label{intro}
Recently the theory of algebras of generalized functions (Colombeau algebras)
\cite{c1,c2} has been restructured to allow for applications in a geometrical
context \cite{found,vim,ndg}. The need for the latter has been clearly demonstrated
by the use of nonlinear generalized function methods in the field of Lie group
analysis of partial differential equations (e.g., \cite{symm, DKP}) and the 
study of singular 
spacetimes in general relativity (see \cite{vickersESI} for a survey).
While diffeomorphism invariance in the so-called full version of Colombeau's
construction (distinguished by a canonical embedding of the space of Schwartz 
distributions) for the scalar case was established in \cite{found,vim} 
using calculus in infinite dimensional (convenient, see \cite{KM}) vector spaces, 
the basic building blocks of the so-called special (or simplified) setting are 
a-priori invariant under the action of a diffeomorphism. The latter, although not 
providing a distinguished embedding of distributions allows to model 
singularities in a 
nonlinear context in a flexible and efficient way. A systematic development
of global analysis using this framework has been started in \cite{RD, ndg, gfvm}.
In the present work we extend this approach by introducing the foundations of a 
generalized pseudo-Riemannian geometry with special emphasis on applications in 
the theory of general relativity. For an alternative approach based on geometry
of vector sheaves (\cite{Mal}), see \cite{MallRos0, MallRos}. 

We start with an in-depth discussion of the notion of invertibility in the algebra of
generalized functions in Section \ref{inv}, thereby setting the stage for the definition
of a generalized pseudo-Riemannian metric on a manifold (Section \ref{genmet}).
In Section \ref{soam} we extend the constructions of \cite{gfvm} 
by introducing the notion 
of a generalized section on a generalized mapping, which, in particular, allows to
define geodesics of a generalized metric.
Section \ref{genprg} is devoted to a systematic investigation of pseudo-Riemannian 
geometry in this setting: in particular we introduce generalized connections and 
prove a ``Fundamental Lemma of \hbox{(pseudo-)} Riemannian geometry'' as well as several 
consistency results with respect to (linear) distributional geometry resp.\ 
the smooth setting. Finally in Section \ref{gengr} after defining the relevant curvature
quantities we give a brief account on applications in general relativity.
The remainder of the present section is devoted to a short review
of Colombeau's construction (in its special variant) and, in particular, the 
global setting introduced in \cite{ndg}. We begin by fixing some notation
from differential geometry.

Throughout this paper $X$ and $Y$ denote paracompact, smooth Hausdorff manifolds 
of dimension $n$ resp.\ $m$. We denote vector bundles with base space $X$ by
$(E,X,\pi_X)$ or simply by $E\to X$ and write a vector bundle chart $(V,\Psi)$
over a chart $(V,\psi)$ of the base $X$ in the form
\beas
 \Psi:\pi^{-1}(V)&\to&\psi(V)\times \K^{n'}\\
  z&\mapsto&(\psi(p),\psi^1(z),\dots,\psi^{n'}(z))\equiv(\psi(p),\bfs\psi(z))\,,
\eeas
where $p=\pi(z)$ and $\K^{n'}$ (with $\K=\R$ or $\K=\C$) is the typical fiber.
Given a vector bundle atlas $(V_\al,\Psi_\al)_\al$ we write the
change of chart in the form $\Psi_\al\circ\Psi_\beta(y,w)=(\psi_{\al\bet}(y),
\bfs{\psi}_{\al\bet}(y)w)$, where $\psi_{\al\bet}:=\psi_\al\circ\psi_\bet^{-1}$
and $\bfs{\psi}_{\al\bet}:\psi_\bet(V_\al\cap V_\bet)\to\mathrm{GL}({n'},\K)$
denotes the transition functions. 

For vector bundles $E\to X$ and $F\to Y$
we denote the space of vector bundle homomorphisms from $E$ to $F$ by 
$\mathrm{Hom}(E,F)$.
Given $f\in\mathrm{Hom}(E,F)$ the induced smooth map on the bases is denoted by
$\underline{f}$, i.e., $\pi_Y\circ f = \underline{f}\circ \pi_X$.
For vector bundle charts $(V,\Phi)$ of $E$ and $(W,\Psi)$ of $F$ we
write the local vector bundle homomorphism
$f_{\Psi\Phi}:=
\Psi\circ f \circ \Phi^{-1}: \vphi(V\cap \underline{f}^{-1}(W))
\times \K^{n'} \to \phi(W) \times \K^{m'}$
in the form
\[
f_{\Psi\Phi}(x,\xi) =
(f_{\Psi\Phi}^{(1)}(x),f_{\Psi\Phi}^{(2)}(x)\cdot\xi)\,.
\]

The space of $\CC^k$-sections of a vector bundle $E\to X$
is denoted by $\Ga^k(X,E)$ and we drop the
superscript in case $k=\infty$. The $(r,s)$-tensor bundle over $X$ will be denoted
by $T^r_s(X)$ and we use the following notation for spaces of tensor
fields  ${\mathcal T}^r_s(X):=\Ga(X,T^r_s(X))$, ${\mathfrak X}:=\Ga(X,TX)$ and
${\mathfrak X}^*:=\Ga(X,T^*X)$, where $TX$ and $T^*X$ denote the tangent
and cotangent bundle of $X$, respectively. For a section $s\in\Ga(X,E)$ we call
$s^i_\al=\Psi^i_\al\circ s\circ \psi_\al^{-1}$ its $i$-th component ($1\leq i\leq {n'}$)
with respect to the vector bundle
chart $(V_\al,\Psi_\al)$. 

The space of $E$-valued distributions of density character
$q$ (see, e.g., \cite{simanca}, Chap.\ 2) 
will be denoted by $\D'(X,E\otimes\Vol^q(X))$ (where $\Vol^q(X)$ is the $q$-volume bundle
of $X$); in particular the space of $(r,s)$-tensor distributions ($q=0$ and $E=T^r_s(X)$) 
will be denoted by $\D'^r_s(X)$.
Whenever convenient we shall use summation convention and abstract index
notation (cf.\ \cite{penrose_rindler}, Chap.\ 2). That is, we denote an $(r,s)$-tensor 
field by $T^{a_1\dots a_r}_{b_1\dots b_s}\in{\mathcal T}^r_s(X)$ while Greek indices, i.e.,
$T^{\al_1\dots \al_r}_{\bet_1\dots \bet_s}$, are used to denote its components with 
respect to a certain basis. Hence equations involving Latin indices are ``true'' 
tensor equations holding in any basis. 

The (special) {\em algebra of generalized functions on $X$} is defined as the
quotient $\gs(X):=\esm(X)/\ns(X)$ of the space $\esm(X)$ of nets of smooth functions 
$(u_\eps)_{\eps\in(0,1]}\in\CC^\infty(X)^{(0,1]}=:{\mathcal E}(X)$ of moderate growth 
modulo the space $\ns(X)$ of negligible nets, where the respective notions of 
moderateness and negligibility are defined (denoting by ${\mathcal P}(X)$ the space of 
linear differential operators on $X$) by the following asymptotic estimates
\beas
        \esm(X)&:=&\{ (u_\eps)_\eps\in{\mathcal E}(X):\ 
        \forall K\subset\subset X,\ \forall P\in{\mathcal P}(X)\ \exists N\in\N:
  \\&&
        \hphantom{mmmmmmmmmmmmmmimmmmm}\sup_{p\in K}|Pu_\eps(p)|=O(\eps^{-N})\}
   \\[1em]
        \ns(X)&:=& \{ (u_\eps)_\eps\in\esm(X):\, 
        \forall K\subset\subset X,\ \forall m \in\N_0:\
        \sup_{p\in K}|u_\eps(p)|=O(\eps^{m}))\}.
\eeas
\noindent
Elements of $\gs(X)$ are denoted by $u=\cl[(u_\eps)_\eps]=
(u_\eps)_\eps+\ns(X)$. $\gs(\_)$ is a {\em fine sheaf of differential algebras} with 
respect to the Lie derivative (w.r.t.\ smooth vector fields) defined by 
$L_\xi u:=\cl[(L_\xi u_\eps)_\eps]$. 

The spaces of moderate resp.\ negligible sequences 
and hence the algebra itself may be characterized locally, i.e., $u\in\gs(X)$ 
if and only if $u_\al:=u\circ\psi^{-1}_\al\in\gs(\psi_\al(V_\al))$ for all charts $(V_\al,\psi_\al)$ 
and $u_\al|_{\psi_\al(V_\al\cap V_
\beta)}=u_\beta|_{\psi_\beta( V_\al\cap V_\beta)}\circ\psi_\beta\circ\psi_\al^{-1}$
for all $\al,\beta$ with $V_\al\cap V_\beta\not=\emptyset$. 

Smooth functions are embedded into $\gs$ simply by the ``constant'' embedding 
$\sigma$, i.e., $\sigma(f):=\cl[(f)_\eps]$, hence $\CC^\infty(X)$ is a  
subalgebra of $\gs(X)$. Moreover, there exist injective sheaf morphisms $\iota:
\D'(\_\,)\hookrightarrow \G(\_\,)$ which coincide with $\sigma$ on $\CC^\infty(\_\,)$
(\cite{ndg} Th.\ 2). These, however, are not canonical (cf.\ the discussion 
in \cite{ndg}, Sec.\ 4). 
In fact, such embeddings depend (in addition to the choice of a mollifier as in the $\R^n$-case)
also on the choice of an atlas for $X$ and of families of cut-off functions, hence are 
non-geometric in an essential way. Nevertheless, there are a number of physically relevant
cases (examples will be given in Section \ref{gengr}) where a distinguished regularization
procedure inducing an embedding into the algebra is suggested by the application under
consideration. Additionally it is often possible to show independence of the distributional
results  results achieved (via the concept of association) through such a procedure from
the choice of embedding. In cases where a canonical embedding is required, on the other hand,
it is possible to employ the intrinsic full Colombeau algebra (providing a canonical embedding
of $\D'$), see \cite{vim}. A theory of generalized pseudo-Riemannian geometry in this
full setting is the subject of ongoing research.

Compatibility with respect to the distributional setting 
is established via the notion of {\em association}, defined as follows:
a generalized function $u$ is called associated with $0$, $u\approx 0$, 
if $\int_X u_\eps \mu \to 0$ ($\eps\to 0$) for all compactly supported one-densities $\mu$ 
and one (hence every) representative $(u_\eps)_\eps$ of $u$. Clearly, $\approx$ 
induces an equivalence relation giving rise to a linear quotient space of $\gs(X)$,
which generalizes the notion of distributional equality to the level of the algebra. 
If $\int_X u_\eps \mu \to \langle w,\mu\rangle$ for some $w\in\D'(X)$ (where $\langle\,,\,\rangle$ 
denotes the distributional action) then $w$ is called the {\em distributional 
shadow} (or {\em macroscopic aspect}) of $u$ and we write $u\approx w$. If it 
exists at all the latter is unique. In the absence of a distinguished embedding it
is useful to define also the stronger notion of {\em $k$-association}. 
We call a generalized function $u$ 
$k$-associated with $0$ ($0\leq k\leq\infty$), $u\approx_k 0$, if for all $l\leq k$,
all $\xi_1,\dots,\xi_l\in {\mathfrak X}(X)$ and one (hence every) representative 
$L_{\xi_1}\dots L_{\xi_l}\,u_\eps\,\to\,0$ uniformly on compact sets.    
Also we say that $u$ admits $f$ as {\em $\CC^k$-associated function},
$u\approx_k f$, if for all $l\leq k$, all $\xi_1,\dots,\xi_l \in {\mathfrak X}(X)$ 
and one (hence any) representative $L_{\xi_1}\dots L_{\xi_l}\,(u_\eps-f)\,\to\,0$
uniformly on compact sets. 
The concept of $k$-association provides a close connection to the classical ($\cc^k$-) 
picture: when modelling singular data in applications it is often possible to 
obtain $\cc^k$-associated functions which serve to establish consistency properties
with the classical setting (cf.\ Proposition \ref{5.4} or Theorem \ref{6.2} below).

Finally, inserting $p\in X$ into $u\in\gs(X)$ yields a well-defined element of the 
ring of constants  ${\mathcal K}$ (corresponding to $\K=\R$ resp.\ $\C$), defined
as the set of moderate nets of numbers ($(r_\eps)_\eps \in \K^{(0,1]}$ with
$|r_\eps| = O(\eps^{-N})$ for some $N$) modulo negligible nets
($|r_\eps| = O(\eps^{m})$ for each $m$). Moreover, generalized functions on $X$ are
characterized by their generalized point values---a feature which distingusihes
them from the purely distributional setting. On $X$ define the space
\begin{equation}
  \label{xtildec}
\widetilde X_c := X_c / \sim  
\end{equation}
of equivalence classes of compactly supported nets $(p_\eps)_\eps
\in X^{(0,1]}$ with respect to the equivalence relation $p_\eps\sim p'_\eps:\Leftrightarrow 
d_h(p_\eps,p'_\eps)=O(\eps^m)$ for all $m$, where $d_h$ denotes the distance 
function on $X$ induced by any Riemannian metric. Then for any generalized 
function $u$ and any $\tilde p \in \widetilde X_c$ the insertion 
$u(\tilde p)$ yields a well defined element of $\tilde{\mathcal K}$ and
$u=0\in\G(X)$ if and only if $u(\tilde p)=0\in{\mathcal K}$ for all generalized 
points $\tilde p\in\widetilde X_c$ (\cite{ndg}, Th.\ 1).

The $\gs(X)$-module of {\em generalized sections} $\Gamma_\gs(X,E)$ of a vector bundle 
$E\to X$ 
is defined along the same lines using analogous asymptotic estimates with 
respect to the 
norm induced by any Riemannian metric on the respective fibers. More precisely,  
setting
$\Gamma_{\mathcal E}(X,E):=(\Ga(X,E))^{(0,1]}$ we define (${\mathcal P}(X,E)$ denoting 
the space of linear differential operators on $\Ga(X,E)$)

\beas
        \Gamma_{\esm}(X,E)&:=& \{ (s_\eps)_{\eps}\in \Gamma_{\mathcal E}(X,E) :
                \ \forall P\in {\mathcal P}(X,E)\, \forall K\comp X \, \exists N\in \N:\\
                 &&\hspace{4cm}\sup_{p\in K}\|Pu_\eps(p)\| = O(\eps^{-N})\}\\
        \Gamma_\ns(X,E)&:=& \{ (s_\eps)_{\eps}\in \Gamma_{\esm}(X,E) :
                \ \forall K\comp X \, \forall m\in \N:\\
                 &&\hspace{3.5cm}\sup_{p\in K}\|u_\eps(p)\| = O(\eps^{m})\}\,,\ \mbox{and finally,}\\
        \Gamma_\gs(X,E)&:=&\Gamma_{\esm}(X,E)/\Gamma_\ns(X,E)\,.
\eeas
\noindent
We denote generalized sections by $s=\cl[(s_\eps)_\eps]=(s_\eps)_\eps+\ns(X,E)$.
Alternatively we may describe a section $s\in\Gamma_\gs(X,E)$ by a family  
$(s_\al)_\al=((s^i_\al)_\al)_{i=1}^{n'}$, 
where $s_\al$ is called the {\em local expression} of $s$ with its {\em components}
$s^i_\al:=\Psi^i_\al\circ s\circ\psi_\al^{-1} \in\gs(\psi_\al(V_\al))$ 
($i=1,\dots, {n'}$) satisfying $s^i_\al(x)\,=\,(\bfs{\psi}_{\al\beta})^i_j
(\psi_\beta\circ\psi^{-1}_\al(x))\,s^j_\beta$ $
(\psi_\beta\circ\psi^{-1}_\al(x))$ for all $x\in \psi_\al(V_\al\cap V_\beta)$, 
where $\bfs{\psi}_{\al\beta}$ denotes the transition functions of the bundle.
Smooth sections of $E\to X$ again may be embedded as constant nets, i.e., 
we define an embedding $\Sigma: \Gamma(X,E) \hookrightarrow \Gamma_\gs(X,E)$ by 
$\Sigma(s)=\cl[(s)_\eps]$. 

$\Ga_\gs(\_\,,E)$ is a fine sheaf of $\gs(X)$-modules. Moreover, the 
$\gs(X)$-module $\gs(X,E)$ is projective and finitely generated (\cite
{ndg}, Th.\ 5). Since $\CC^\infty(X)$ is a subring of $\gs(X)$,
$\Gamma_\gs(X,E)$ also may be viewed as $\CC^\infty(X)$-module and the two respective module 
structures are compatible with respect to the embeddings. Furthermore we have 
the following algebraic characterization of the space of generalized sections (\cite{ndg},
Th.\ 4) 
\begin{equation}\label{tensorp}
  \Gamma_\gs(X,E)=\gs(X)\otimes\Ga(X,E)\,,
\end{equation}
where the tensor product is taken over the module $\CC^\infty(X)$. 

Compatibility with respect to the classical resp.\ distributional
setting again is accomplished using the concept of ($k$-)association. 
A section $s\in\Gamma_\gs(X,E)$ is called {\em associated} with $0$, $s\approx 0$, if 
all its components $s^i_\al\approx 0$ in $\gs(\psi_\al(V_\al))$. $s$ allows for 
the {\em distributional shadow} $w$, $s\approx w\in\D'(X,E)$, if $s^i_\al\approx w^i_\al\in
\gs(\psi_\al(V_\al))$ for all $\al,i$ (and $w^i_\al$ denoting the local expression
of the distribution $w$). 
Similarly  $s$ is called {\em $\CC^k$-associated} with $0$ $(0\leq k\leq\infty)$, $s\approx_k 0$, if for 
one (hence every) representative $(s_\eps)_\eps$ all components $s^i_\al\,_\eps\to 0$ 
uniformly on compact sets in all derivatives of order less or equal to $k$.
We say that $s$ allows $t\in\Ga^k(X,E)$ as a {\em $\CC^k$-associated section},
$s\approx_k t$, if  for  one (hence every) representative $(s_\eps)_\eps$  all components
$s^i_\al\,_\eps\to t^i_\al$ uniformly on compact sets in all derivatives of order 
less or equal to $k$.

Generalized tensor fields (i.e., elements of $\gs^r_s(X):=\Ga_\gs(X,T^r_s(X))$) 
may be viewed likewise as $\CC^\infty(X)$-multilinear mappings taking smooth 
vector fields resp.\ one-forms to $\gs(X)$ or as $\gs(X)$-multilinear mappings taking 
generalized vector resp.\ covector fields to generalized functions, i.e., as 
$\CC^\infty(X)$- resp.\ $\gs(X)$-modules we have (\cite{ndg}, Th.\ 6)
\beas
   \gs^r_s(X)&\cong&L_{\CC^\infty(X)}({\mathfrak X}^*(X)^r,{\mathfrak X}(X)^s;\gs(X))\\
   \gs^r_s(X)&\cong&L_{\gs(X)}(\gs^0_1(X)^r,\gs^1_0(X)^s;\gs(X)).
\eeas
Given a generalized tensor field $T\in{\gs}^r_s(X)$ 
we shall call the $n^{r+s}$ generalized functions on $V_\al$ defined by
\[ 
        T^\alpha\,^{\beta_1\dots \beta_r}_{\gamma_1\dots \gamma_s}
        \,:=\,T|_{V_\al}(dx^{\beta_1},\dots ,dx^{\beta_r},\pa_{\gamma_1},\dots,\pa_{\gamma_s})
\]
its {\em components} with respect to the chart $(V_\al,\psi_\al)$.

In \cite{ndg} many concepts of classical differential geometry (in particular,
Lie derivatives with respect to both smooth and generalized vector fields, 
Lie brackets, 
tensor product, contraction, exterior algebra etc.) have been generalized to this new 
setting and will be used in the sequel. 

\section{Invertibility}\label{inv}
Prior to our analysis of generalized semi-Riemannian geometry, in 
the present section we are going to derive several 
characterization results concerning invertibility in the Colombeau algebra
which will be essential for the algebraic aspects of the theory to be developed
in the subsequent sections.

We begin with a characterization of multiplicative invertibility in $\gs(X)$.
\bp \label{invchar}
Let $u\in \gs(X)$. The following are equivalent:
\begin{itemize}
\item[(i)] There exists $v\in \gs(X)$ with $uv = 1$.
\item[(ii)] For each representative $(u_\eps)_\eps$ of $u$ and each $K\comp X$
there exist $\eps_0>0$ and $q\in \N$ such that $\inf_{p\in K} |u_\eps(p)| \ge \eps^q$
for all $\eps<\eps_0$.
\end{itemize}
\end{proposition}

\pr (i) $\Rightarrow$ (ii): Let $u = \cl[(u_\eps)_\eps]$ and $v=\cl[(v_\eps)_\eps]$.
By assumption there exists $N = \cl[(n_\eps)_\eps] \in \ns(X)$ 
such that $u_\eps v_\eps \equiv 1 + n_\eps$.
We first claim that there exists some $\eps_0$ such that $v_\eps(p)\not=0$ for all
$p\in K$ and all $\eps < \eps_0$. Indeed, otherwise there would exist a zero-sequence
$\eps_m$ and a sequence $p_m$ in $K$ such that $v_{\eps_m}(p_m) = 0$ for all $m$.
But then $0 = u_{\eps_m}(p_m)v_{\eps_m}(p_m) = 1 + n_{\eps_m}(p_m) \to 1 \ 
(m\to \infty)$, a contradiction. Since $(v_\eps)_\eps \in \esm(X)$ there exists 
$l\in \N$ and $\eps_1>0$ with $\sup_{p\in K} |(v_\eps)_\eps(p)| < \eps^{-l}$ for 
$\eps<\eps_1$. Hence
$$
|u_\eps(p)|> \eps^l(1-|n_\eps(p)|) > \eps^{l+1}
$$ 
for $\eps$ small, uniformly for $p\in K$.

(ii) $\Rightarrow$ (i): Let $(X_m)_{m\in\N}$ be an open covering of $X$ consisting
of relatively compact sets. Then by assumption, $v_\eps^m(p) := \frac{1}{u_\eps(p)}$
exists for $p\in X_m$ and $\eps$ sufficiently small (for all other $\eps$
we may set $v_\eps^m\equiv 0$). Also, $(v_\eps^m)_\eps \in \esm(X_m)$, so 
$v^m := \cl[(v_\eps^m)_\eps]$ is a well-defined element of $\gs(X_m)$. By definition,
for $X_m\cap X_k \not=\emptyset$ we have $v^m|_{X_m\cap X_k} = v^k|_{X_m\cap X_k}$,
so $\{v^m|m\in \N\}$ forms a coherent family. Since $\gs(\_)$ is a sheaf of differential
algebras (cf.\ \cite{RD}, \cite{ndg}) it follows that
there exists a unique element $v\in\gs(X)$ with $v|_{X_m} = v_m$
for all $m\in \N$. It is clear that $v$ is the desired multiplicative inverse of $u$.  
\ep

We call an element 
$r\in \gK$ {\em strictly nonzero} 
if there exists some representative $(r_\eps)_\eps$
and an $m\in \N$ with $|r_\eps|\ge \eps^m$ for $\eps$ sufficiently small.
By specializing \ref{invchar} to $K=\{p\}$ ($p\in X$) it follows that 
$r$ is invertible if and only if it is strictly nonzero.
Also, \ref{invchar} can be restated as follows: a Colombeau function possesses
a multiplicative inverse if and only if it is strictly nonzero, uniformly on compact sets.
\bp \label{noninvertible}
Let $r\in \gK$. The following are equivalent:
\begin{itemize}
\item[(i)] $r$ is not invertible.
\item[(ii)] $r$ is a zero divisor. 
\end{itemize}
\end{proposition}

\pr
(i) $\Rightarrow$ (ii): By the above we have: $\forall N\in \N$ $\forall \eps_0$
$\exists \eps<\eps_0$: $|r_\eps| < \eps^N$. Thus there exists a sequence $\eps_k \searrow
0$ with $|r_{\eps_k}| < \eps_k^k$. Define $(a_\eps)_\eps$ by $a_\eps = r_{\eps_{k+1}}$
for $\eps_{k+1} \le \eps <\eps_k$. Then
$$
|a_\eps| = |r_{\eps_{k+1}}| < \eps_{k+1}^{k+1} \le \eps^{k+1} \qquad 
(\eps_{k+1} \le \eps <\eps_k)\,.
$$
Hence $|a_\eps| < \eps^{N+1}$ for $0<\eps <\eps_N$, i.e., $(a_\eps)_\eps \in \ns$.
Setting $\tilde r_\eps = r_\eps - a_\eps$ we obtain a representative $\tilde r$ of
$r$ with $\tilde r_{\eps_k}=0$ for all $k\in \N$.
Denote by $s$ the class of $(s_\eps)_\eps$, where
$$
s_\eps = \left\{
\begin{array}{l}
1, \quad \eps \in \{\eps_1,\eps_2,\dots\} \\
0, \quad  \eps \not\in \{\eps_1,\eps_2,\dots\} 
\end{array}
\right.
$$
Then $\tilde r_\eps s_\eps \equiv 0$, so $r s = 0$, but $s\not=0$ in $\gK$.\\
(ii) $\Rightarrow$ (i) is obvious.
\ep

A comprehensive study of algebraic properties of $\gK$ can be found in \cite{AJ}.
Clearly the analogue of \ref{noninvertible} is wrong for ${\mathcal C}^\infty(X)$. 
It is also false for $\gs(X)$:

\bex \label{zdex} {\rm 
Let $X=\R$ and let $x$ denote the identical function on $\R$, considered as an
element of $\gs(\R)$. Then $x$ is not invertible by \ref{invchar} and we are going to
show that it is not a zero divisor in $\gs(\R)$. To this end, let $u= 
\cl[(u_\eps)_\eps]\in \gs(\R)$ such that $xu = 0$ in $\gs(\R)$. Suppose that 
$(u_\eps)_\eps \not\in \ns(\R)$. Since $(u_\eps)_\eps$ necessarily satisfies
the $\ns$-estimates on every compact set not containing $0$ it follows that
there exists $K\comp \R$ with $0\in K$, $q_0\in \N$, $\eps_m \searrow 0$ and $x_m\in K$
such that $|u_{\eps_m}(x_m)| > \eps_m^{q_0}$ for all $m\in \N$. 
Without loss of generality we may suppose that $x_m\to 0$. For any $r\in \N$ and 
$m$ large we have $|x_m u_{\eps_m}(x_m)|\le\eps_m^r$. Hence $|x_m| 
\le \eps_m^{r-q_0}$ for these $m$, i.e., $x_m$ goes to $0$ faster than any power
of $\eps$. Now
$$
|u_{\eps_m}(x_m)| \le |x_m u_{\eps_m}'(x_m)| + |(xu_{\eps_m})'(x_m)| 
$$
Since $(u_\eps)_\eps$ is moderate there exists $N\in \N$ such 
that $|u_{\eps_m}'(x_m)|\le \eps_m^{-N}$. By the above, $|x_m|\le \eps_m^{N+q_0+1}$
for $m$ large and so, taking into account the $\ns$-estimate for
$(xu_\eps)'$ we obtain
$$
|u_{\eps_m}(x_m)| \le \eps_m^{q_0+1}   \qquad (m \ \mbox{large})
$$
Consequently, $\eps_m^{q_0} \le \eps_m^{q_0+1}$ for large $m$, a contradiction.}
\end{example}

We note that the notion of zero divisor in ${\mathcal C}^\infty(X)$ and $\gs(X)$
differs: in fact (with $X=\R$) $e^{-1/x^2}$ is obviously not a zero 
divisor in ${\mathcal C}^\infty(\R)$ but $e^{-1/x^2} \cdot \delta(x) = 0$ in $\gs(\R)$,
so it is a zero divisor in $\G(X)$.

An element $f$ of ${\mathcal C}^\infty(X)$ is invertible if and only if $f(p)$ is invertible
in $\R$ for each $p\in X$. Our next aim is to find the appropriate generalization
of this observation to the context of Colombeau algebras. A straightforward adaptation
of the smooth case is impossible as is demonstrated by the following example:
\bex {\rm Let $X=\R$ and set 
$$
u_\eps(x) := \eps^{\frac{x^2}{x^4+\eps^4}}\,.
$$
The net $(u_\eps)_\eps$ is moderate since $|u_\eps(x)| \le 1$ and  
$u_\eps^{(k)}(x) = r_k(x,\eps,\ln(\eps))u_\eps(x)$ with $r_k$ a rational
function ($k\in \N$). Thus $u:= \cl[(u_\eps)_\eps]$ is a well-defined element
of $\gs(\R)$. We are going to show that $u(x_0)$ is invertible in $\gK$ for each
$x_0\in \R$ but that $u$ is not invertible in $\gs(\R)$.

In fact, $u_\eps(0)=1$ for all $\eps$ and for $x_0\not=0$ we have 
$u_\eps(x_0) \ge \eps^{\frac{1}{x_0^2}}$, so $u(x_0)$ is 
strictly nonzero, hence invertible in $\gK$ for each $x_0 \in \R$. However,
$u$ is not {\em uniformly} strictly nonzero on compact sets, hence not invertible
in $\gs(\R)$ by \ref{invchar}. To see this, take $K=[0,1]$. Then $x=\eps \in K$ but
$u_\eps(\eps) = \eps^{\frac{1}{2\eps^2}} \ge \eps^q$ cannot hold
identically in $\eps$ in any interval $(0,\eps_0)$ for any $q\in \N$.}
\end{example}
 
The correct generalization of the result in the smooth case uses the point value
characterization of elements of $\gs(X)$ derived in \cite{point}, \cite{ndg}:
\bp \label{invertprop} 
$u\in \gs(X)$ is invertible if and only if $u(\tilde p)$ is invertible in $\gK$ for
each $\tilde p \in \widetilde X_c$.
\end{proposition}

\pr The condition is obviously necessary. Conversely, suppose that $u\in \gs(X)$ 
is not invertible. Then by \ref{invchar} there exist $K\comp X$ and sequences
$\eps_m\to 0$, $q_m \to \infty$ and $p_m \in K$ such that
$$
|u_{\eps_m}(p_m)| < \eps_m^{q_m} \qquad (m\in \N)\,.
$$
Set $p_\eps := p_m$ for $\frac{1}{m+1} < \eps \le \frac{1}{m}$. Then $\tilde p :=
\cl[(p_\eps)_\eps]$ defines a compactly supported generalized point in $X$, i.e.,
an element of $\widetilde X_c$. 
By construction, $u(\tilde p)$ is not invertible in $\gK$.
\ep

Finally, we shall need the following characterization of nondegeneracy in $\gK^n$:
\blem \label{nondeglemma}
 Let $A \in \gK^{n^2}$. The following are equivalent:
\begin{itemize}
\item[(i)] $A$ is nondegenerate, i.e., $\xi\in \gK^n$, $\xi^tA\eta = 0$ $\forall
\eta \in \gK^n$ implies $\xi=0$.
\item[(ii)] $A: \gK^n \to \gK^n$ is injective.
\item[(iii)] $A: \gK^n \to \gK^n$ is bijective.
\item[(iv)] $\det(A)$ is invertible.
\end{itemize}
\end{lemma}

\pr
(ii) $\Leftrightarrow$ (iii) $\Leftrightarrow$ (iv): 
By \cite{bourbaki-algebra},
Ch.\ III, \S 8, Prop.\ 3, (ii) is equivalent with $\det(A)$ 
not being a zero divisor in $\gK$. 
\cite{bourbaki-algebra}, Ch.\ III, \S 8, Th.\ 1 shows that 
(iii) is equivalent with $\det(A)$ being invertible.
Hence the claim follows from \ref{noninvertible}.\\
(i) $\Leftrightarrow$ (ii): First (i) is equivalent with $A^t$ being
injective. Indeed $\xi^tA\eta = 0\ \forall\eta \in \gK^n\ \Leftrightarrow\ 
\xi^t A = 0$ (just set $\eta = e_i$, the $i$-th unit vector in $\K^n
\subseteq \gK^n$), which in turn is equivalent with $A^t\xi = 0$. 
Since $\det(A) = \det(A^t)$ the claim follows from (ii) 
$\Leftrightarrow$ (iv).
\ep

For further studies of linear algebra over the ring $\gK$ we refer to \cite{LT}.

\section{Generalized Metrics}\label{genmet}

In \cite{ndg}, Th.\ 7 the following isomorphism of $\G(X)$-modules was established:
\begin{equation} \label{gsrmodule}
\gs^r_s(X) \cong L_{\gs(X)}(\gs^0_1(X)^r,\gs^1_0(X)^s;\gs(X))
\end{equation}
We will make use
of this identification in the following characterization result which will motivate 
our definition of generalized metrics.
\bt \label{mainmetric}
Let $\hat g\in \gs^0_2(X)$. The following are equivalent:
\begin{itemize}
\item[(i)] For each chart $(V_\al,\psi_\al)$ and each $\tilde x \in 
(\psi_\al(V_\al))^{\sim}_c$ (cf.\ (\ref{xtildec})) the map
$\hat g_\al(\tilde x): \gK^n\times\gK^n\to\gK^n$ is symmetric and nondegenerate.
\item[(ii)] $\hat g:
\gs^1_0(X) \times \gs^1_0(X) \to \G(X)$ is symmetric and
$\mbox{det}(\hat g)$ is invertible in $\G(X)$.
\item[(iii)] $\mbox{det}(\hat g)$ is invertible in $\G(X)$ and for each relatively
compact open set $V\subseteq X$ there exists a representative $(\hat g_\eps)_\eps$ of
$\hat g$ and an $\eps_0 > 0$ such that $\hat g_\eps|_V$ is a smooth pseudo-Riemannian
metric for all $\eps<\eps_0$.
\end{itemize}
\end{theorem}

\pr 
(i) $\Leftrightarrow$ (ii): Supposing (i), for any $\xi, \eta \in \G^1_0(\psi_\al(V_\al))$ 
it follows from \cite{point}, Th.\ 2.4\ and \cite{gfvm}, Prop.\ 3.8 that\ 
$\hat g_\al(\xi,\eta)\,=\,
\hat g_\al(\eta,\xi)$,\ so\ $\hat g|_{V_\al} \in$\ $ L_{\gs(V_\al)}\ 
(\gs^1_0(V_\al),\gs^1_0(V_\al);\\ \gs(V_\al))$ is symmetric. Since 
$L_{\gs(\_)}$$(\gs^0_1(\_)^r,$$\gs^1_0(\_)^s;$ $\gs(\_))$ is a sheaf, 
sym\-metry of 
$\hat g$ follows. Moreover, by \ref{nondeglemma} and \ref{invertprop} it follows
that $\det(\hat g)|_{V_\al}$ is invertible for each $\al$, so $\det(\hat g)$ is
invertible on $X$.\\
The converse direction follows immediately from \ref{invertprop}  and \ref{nondeglemma}.\\
(ii) $\Rightarrow$ (iii):
We first note that for $\xi,\eta \in \gs^1_0(X)$, 
$$
\hat g(\xi,\eta) = \frac{1}{2}(\hat g(\xi,\eta) + \hat g(\eta,\xi)) + 
\frac{1}{2}(\hat g(\xi,\eta) - \hat g(\eta,\xi))\,.
$$
The second term in this expression is $0$ in $\gs(X)$ for all $\xi,\eta$ by assumption.
Hence $\hat g$ equals the element of $\gs^0_2(X)$ corresponding via
(\ref{gsrmodule}) to $(\xi,\eta)\to \frac{1}{2}$ $(\hat g(\xi,\eta)$ $+\hat g(\eta,\xi))$. 
From this we obtain a representative $(\hat g_\eps)_\eps$ of $\hat g$ such that each 
$\hat g_\eps: \mathfrak{X}(X)\times \mathfrak{X}(X) \to {\mathcal C}^\infty(X)$ is symmetric. 
Moreover, by \ref{invchar} for any $K\subset\subset X$ there exists $\eps_0 > 0$, $q\in \N$ 
such that $\inf_{p\in K} |\det(\hat g_\eps(p))| > \eps^q$ for $\eps<\eps_0$. In particular,
each $\hat g_\eps$ is nondegenerate, hence a pseudo-Riemannian metric on 
any open $V\subset K$.\\
(iii) $\Rightarrow$ (i): Let $\tilde x \in \psi_\al(V_\al)^{\sim}_c$ be supported in
$K \comp \psi_\al(V_\al)$ 
and choose a representative $(\hat g_\eps)_\eps$ of $\hat g$ such that
each $\hat g_\eps$ is a pseudo-Riemannian metric on a neighborhood of $\psi_\al^{-1}(K)$. 
Since each $\hat g_\eps$ is symmetric it follows that 
$\hat g_\al(\tilde x) = \cl[(\hat g_{\al\eps}(x_\eps))_\eps]$
is symmetric as a map from $\gK^n\times\gK^n$ to $\gK$. 
Finally, nondegeneracy follows from \ref{nondeglemma}.
\ep

\bd\label{gmetric}
Suppose that $\hat g\in{\gs}^0_2(X)$ satisfies one (hence all) of the 
equivalent conditions in \ref{mainmetric}.
If there exists some $j\in \N_0$ such that for each relatively compact open set $V
\subseteq X$ there exists a representative $(\hat g_\eps)_\eps$ of $\hat g$ as in 
\ref{mainmetric} (iii) such that the index of each $\hat g_\eps$ equals $j$ we call $j$ the
{\em index} of $\hat g$.
\end{definition}

For the above notion of index to make sense we have to establish that it does not
depend on the representative of $\hat g$ used in (iii). To secure this property
we make use of a result from perturbation theory of finite dimensional linear
operators.

\bp The index of $\hat g\in{\gs}^0_2(X)$ as introduced above is well-defined.
\end{proposition}

\pr Let $V\subseteq X$ be relatively compact, let $(\hat g_\eps)_\eps$ 
be a representative 
of $\hat g$ as above and denote by $\hat\la^1_\eps\geq\dots\geq\hat\la^n_\eps$ its 
eigenvalues. By
\ref{mainmetric} (iii) and \ref{noninvertible} it follows that 
each $\hat \la^i$ is invertible in $\gK$ (otherwise $\det\hat g$ would be a zero 
divisor).  Hence each $\hat \la^i$ is strictly nonzero, i.e., 
there exists $r\in\N_0$ such that 
\beq\label{bdafz}
|\hat\la^i_\eps|>\eps^r
\eeq
for all $1\leq i\leq n$ and $\eps$ small. 
Let $(\tilde g_\eps)_\eps$ be another representative of 
$\hat g$ as in \ref{mainmetric} (iii) with eigenvalues 
$\tilde\la^1_\eps\geq\dots\geq\tilde\la^n_\eps$. 
By \cite{bhatia}, 3, Th.\ 8.1 we have $\max_i|\tilde\la^i_\eps-\hat\la^i_\eps|\leq||\tilde g_\eps-
\hat g_\eps||$. Hence from the $\ns$-estimate we conclude that for all $i$ $(1\leq i\leq n)$, 
$|\tilde\la^i_\eps-\hat\la^i_\eps| = O(\eps^m)$ $\forall m$. Thus by (\ref{bdafz}) $\tilde\la^i_\eps$
and $\hat\la^i_\eps$ have the same sign for small $\eps$.
\ep

\bd \label{gmdef} \
\begin{itemize}
\item[(i)] 
 A generalized $(0,2)$-tensor field $\hat g\in{\gs}^0_2(X)$
\hspace*{-3pt} is called a 
{\em generalized  \hbox{(pseudo-)} Riemannian metric} if it 
satisfies one of the equivalent conditions in \ref{mainmetric} and 
possesses an index (different from $0$).
\item[(ii)] 
We call a paracompact, smooth Hausdorff manifold
$X$ furnished with a generalized (pseudo-)Rie\-mann\-ian metric $\hat g$ a
{\em generalized (pseudo-)Rieman\-nian manifold} or, if the index 
of $\hat g$ is $1$ or 
$n-1$, a {\em generalized spacetime} and denote it by $(X,\hat g)$. 
The action of the metric on a pair of generalized vector fields will
be denoted by $\hat g(\xi,\eta)$ and $\langle\xi,\eta\rangle$, equivalently. 
\end{itemize}
\end{definition}

\brem {\rm
Let us compare this notion of generalized metric with the ones
introduced in the purely distributional picture in 
\cite{marsden}, 10.6, and in \cite{parker}. 
In \cite{marsden}, a distributional $(0,2)$-tensor field
$g\in{\D'}^0_2(X)$ is called nondegenerate if $g(\xi,\eta) = 0$ for all
$\eta\in \mathfrak{X}(X)$ implies $\xi=0\in\mathfrak{X}(X)$. 
However, this ``nonlocal'' condition is too weak to reproduce
the classical notion of nondegeneracy; just take $ds^2 = x^2\,dx^2$.
In \cite{parker}, on the other hand, $g\in
{\D'}^0_2(X)$ is called nondegenerate if it is nondegenerate 
(in the classical sense) off its singular support, so in this approach 
no statement at all is made at the singularities of the metric. 

Even by combining both notions, i.e., by calling a distributional 
$(0,2)$-tensor field nondegenerate if it satisfies both conditions one 
arrives at a comparatively weak notion. 
To see this, take $ds^2 = (x^2 + \delta(x))\, dx^2$
on $\R$. This metric is easily seen to be nondegenerate in the above sense.
According to  \ref{gmdef}, however, it depends on the ``microstructure''
of $\delta$, i.e., on the chosen regularization of $\delta$ whether 
or not $ds^2$ is nondegenerate in the $\gs$-setting. }
\end{remark}

By \ref{invchar} (i), for any representative $(\hat g_\eps)_\eps$ of a generalized
pseudo-Riemannian metric we have
\begin{equation}
  \label{invmetform}
\forall K \subset\subset X \ \exists m\in\N \ \exists \eps_0>0: 
\inf\limits_{p\in K}|\det \hat g_\eps(p)|\geq\eps^m \ \forall \eps<\eps_0.
\end{equation}

\bp Let $(X,\hat g)$ be a generalized pseudo-Riemannian manifold and let 
$\hat g = \cl[(\hat g_\eps)_\eps]$. Then the {\em inverse metric}
$\hat g^{-1}:= \cl[(\hat g^{-1}_\eps)_\eps]$ is a well-defined 
element of ${\gs}^2_0(X)$.
\end{proposition}

\pr 
We first note that by (\ref{invmetform}), in any  relatively compact
chart $V$ and for fixed $\eps < \eps_0(V)$ we 
may define (in usual notation) $\hat g^{ij}\,_\eps$ to be the pointwise
inverse of $\hat g_{ij}\,_\eps$ which is obviously a smooth 
$(2,0)$-tensor field.  
By the cofactor formula of matrix inversion we have 
$\hat g^{ij}\,_\eps=\mathrm{cof}(\hat g_{ij}\,_\eps)/(\det\hat g_{ij}\,_\eps)$,
so from (\ref{invmetform}) we conclude that $(\hat g^{ji}\,_\eps)_\eps$  is moderate
on its domain of definition. By the same reasoning, choosing a different 
representative of $\hat g$ in this calculation perturbs $\hat g^{ij}$ merely by an
element of $\ns$. From the sheaf property of $\gs^2_0(\_)$ the 
result follows.
\ep

From now on we denote the inverse metric (using abstract index notation, 
cf.\ Section \ref{intro}) by $\hat g^{ab}$, its components by $\hat g^{ij}$ 
and the components of a representative by $\hat g^{ij}\,_\eps$. Also, we shall 
denote the 
line element 
by $\hat{ds}^2=\cl[(\hat{ds}_\eps^2)_\eps]$.

\bex\label{mex}\ 
{\rm
\begin{itemize}
\item[(i)] 
A sufficient condition for a sequence $(g_\eps)_\eps$ 
of classical (smooth) metrics of constant index to
constitute a representative of a generalized metric---apart from being moderate---is 
to be zero-associated (i.e., to converge locally uniformly) to a classical 
(then necessarily continuous) metric $g$. Indeed, (\ref{invmetform}) is satisfied
in this case since $\det g_\eps \to \det g$ uniformly on compact sets, so the
claim follows from \ref{mainmetric} (iii). 
\item[(ii)] 
The metric of a two-dimensional cone was modelled in \cite{clarke} by a 
generalized metric (in the full setting) obtained through embedding via
convolution. 
\item[(iii)]
The
line element of impulsive pp-waves in \cite{geo,geo2,penrose} was modelled by
\begin{equation}
   \hat{ds}^2=f(x,y)D(u)du^2-dudv+dx^2+dy^2\,,
\end{equation}
where $D$ denotes a generalized delta function which allows for a strict delta
net as a representative (see Section \ref{gengr}) and $f$ is a smooth function.
\item[(iv)] Further examples may be found e.g., in \cite{hb5,ultra,MaN}.
\end{itemize}
}
\end{example}

Since taking the determinant is a polynomial operation we cannot
expect association to be compatible with inverting a metric. However, the
analogous statement for $k$-association holds by an application of \cite{ndg},
Prop.\ 3 (ii).

\bp
Let $\hat g_{ab}$ a generalized metric and $\hat g_{ab}\approx_k g_{ab}$, where
$g_{ab}$ is a classical $\CC^k$-pseudo-Riemannian metric. Then
$\hat g^{ab}\approx_k g^{ab}$.
\end{proposition}

Additional important properties of generalized metrics are presented in
the following result.
\bp \label{mlemma}
Let $(X,\hat g)$ be a generalized pseudo-Riemannian manifold.
\begin{itemize}
\item[(i)] 
$\xi\in \gs^1_0(X)$, $\hat g(\xi,\eta)=0$ $\forall \eta\in{\gs}^1_0(X)$ 
$\Rightarrow$ $\xi=0$.
\item[(ii)] $\hat g$ induces a $\gs(X)$-linear isomorphism 
${\gs}^1_0(X)\to{\gs}^0_1(X)$
by $\xi\mapsto \hat g(\xi,\,.\,)$.
\end{itemize}
\end{proposition}

\pr (i) We have to show that for $\hat g=\cl[(\hat g_\eps)_\eps]$, 
$\xi=\cl[(\xi_\eps)_\eps]$, $(\hat g_\eps(\xi_\eps,\eta_\eps))_\eps\in\ns(X)$ 
for all $\eta=\cl[(\eta_\eps)_\eps]$ 
implies $(\xi_\eps)_\eps\in({\ns})^1_0(X)$.
By the sheaf properties of $\gs^r_s(\_)$ it suffices to establish the claim 
locally on every chart $V_\al$, i.e.,
\[ 
        \left(\hat g_{ij}\,_\eps\xi^i_\eps\eta^j_\eps\right)_\eps
        \in\ns(V_\al)\ \forall (\eta^j_\eps)_\eps
        \,\Rightarrow\, (\xi^i\,_\eps)_\eps\in\ns(V_\al)\quad i=1\dots n\,.
\]
Setting $\eta^j\,_\eps=\sum_l\hat g^{jl}\,_\eps\,\xi^l_\eps$ gives $(\sum_{ij}\hat g_{ij}\,_\eps
\xi^i_\eps\eta^j_\eps)=\sum_i(\xi^i\,_\eps)^2$. 
It follows that each $\xi^i$ satisfies the ${\mathcal N}$-estimates of order $0$. Thus
the claim follows from \cite{found}, Th.\ 13.1.\\
(ii) By (\ref{gsrmodule}) $\xi^*:=\hat g(\xi,.)$ is indeed a 
one-form and the assignment $\xi\mapsto\xi^*$ is $\gs(X)$-linear.
Moreover, injectivity of this map follows from \ref{mlemma} (i). 
It remains to show that
the assignment is onto. Locally, any generalized one-form can be written as 
$A =A_idx^i$. Define a
generalized vector field by $V=\hat g^{jk}\,A_j\,\pa_k$. Then
\[ 
        \langle V,\pa_l\rangle
        \,=\,\hat g^{jk}\,A_k\,\langle\pa_j,\pa_l\rangle
        \,=\,A_l
\] 
and the result again follows from the sheaf property of $\gs^1_0(\_)$.
\ep

The isomorphism in (ii) above---as in the classical context---extends naturally to 
generalized tensor fields of higher types. Hence from now on we shall use the 
common conventions on upper and lower indices also
in the context of generalized tensor fields. In particular, identifying
a vector field $\xi^a\in{\gs}^1_0(X) $ with its metrically equivalent one-form 
$\xi_a$ we denote its contravariant respectively covariant components
by $\xi^i$ and $\xi_i$. A similar convention will apply to representatives.

\section{Sections on a generalized mapping}\label{soam}
In \cite{gfvm}, the space $\gs[X,Y]$ of Colombeau generalized functions on
the manifold $X$ taking values in the manifold $Y$ as well as the space 
$\mathrm{Hom}_{\gs}[E,F]$
of generalized vector bundle homomorphism from $E$ to $F$ was defined. 
In order to obtain a consistent description of geodesics of generalized 
pseudo-Riemannian metrics we need some additional constructions extending the
framework introduced there. The present section is devoted to
the development of these concepts and we shall freely use notations and definitions
from \cite{gfvm}. For the convenience of the reader, however, we recall the definition
of generalized functions valued in a manifold.
\bd Let $X$, $Y$ be paracompact, smooth Hausdorff manifolds.
\begin{itemize}
\item [(a)] The space $\esm[X,Y]$ of compactly bounded
(c-bounded, for short) moderate maps from $X$ to $Y$ is defined as the set
of all $(u_\eps)_\eps \in \cinfty(X,Y)^{(0,1]}$ such that
\begin{itemize}
  \item[(i)] $\forall K\comp\Om\ \exists \eps_0>0\  \exists K'\comp Y \  \forall
             \eps<\eps_0:\ u_\eps(K) \subseteq K'$.
  \item[(ii)]
   $\forall k\in\N$,
   for each chart
   $(V,\vphi)$
   in $X$, each chart $(W,\psi)$ in $Y$, each $L\comp V$
   and each $L'\comp W$
   there exists $N\in \N$  with
   $$\sup\limits_{p\in L\cap u_\eps^{-1}(L')} \|D^{(k)}
   (\psi\circ u_\eps \circ \vphi^{-1})(\vphi(p))\| =O(\eps^{-N}).
   $$
\end{itemize} 
\item[(b)] Two elements $(u_\eps)_\eps$, $(v_\eps)_\eps$ of $\esm[X,Y]$ are
called equivalent, $(u_\eps)_\eps \sim (v_\eps)_\eps$, 
if the following conditions are satisfied: 
\begin{itemize}
\item[(i)]  For all $K\comp X$, $\sup_{p\in K}d_h(u_\eps(p),v_\eps(p)) \to 0$
($\eps\to 0$)
for some (hence every) Riemannian metric $h$ on $Y$.
\item[(ii)] $\forall k\in \N_0\ \forall m\in \N$,
for each chart
   $(V,\vphi)$
   in $X$, each chart $(W,\psi)$ in $Y$, each $L\comp V$
   and each $L'\comp W$:
$$
\sup\limits_{p\in L\cap u_\eps^{-1}(L')\cap v_\eps^{-1}(L')}\!\!\!\!\!\!\!\!\!\!\!\!\!\!\!\!\!\!\!
\|D^{(k)}(\psi\circ u_\eps\circ \vphi^{-1}
- \psi\circ v_\eps\circ \vphi^{-1})(\vphi(p))\|
=O(\eps^m).
$$
\end{itemize}
\item[(c)] The quotient $\gs[X,Y]:=\esm[X,Y]/\sim$ is called the space
of Colombeau generalized functions on $X$ {\em valued in} $Y$. 
\end{itemize}
\end{definition}

The space $\mathrm{Hom}_{\gs}[E,F]$ of generalized vector bundle 
homomorphisms is defined along the same lines (see \cite{gfvm}, Section 3). 
Here we introduce a ``hybrid'' variant of generalized mappings 
defined on a manifold and taking values in a vector bundle. 
\bd \label{hybridmod} 
Let $(F,\pi_Y,Y)$ be a vector bundle and denote by $\esmh[X,F]$ the set of
all nets $(u_\eps)_\eps \in \cinfty(X,F)^{(0,1]}$ satisfying (with 
$\underline{u_\eps} := \pi_Y\circ u_\eps$)
\begin{itemize}
\item[(i)] 
$\forall K\comp X\ \exists K'\comp Y\ \exists \eps_0>0\ \forall \eps<\eps_0$
$\underline{u_\eps}(K) \subseteq K'$.
\item[(ii)] $\forall k\in \N_0
\ \forall (V,\vphi)$ chart in $X$ $\forall (W,\Psi)$ vector bundle
chart in $F$ $\forall L\comp V\ \forall L'\comp W
\ \exists N\in \N\ \exists \eps_1>0\ \exists C>0$ such that
$$
\|D^{(k)}(\Psi\circ u_\eps \circ \vphi^{-1})(\vphi(p))\| \le
C\eps^{-N}
$$
for each $\eps<\eps_1$ and each $p\in
L\cap\underline{u_\eps}^{-1}(L')$.
\end{itemize}
\end{definition}

In particular, $(u_\eps)_\eps \in \esmh[X,F]$ implies $(\underline{u_\eps})_\eps
\in \esm[X,Y]$.
\bd \label{hybridequiv} 
$(u_\eps)_\eps$, $(v_\eps)_\eps \in \esmh[X,F]$ are called equivalent,
$(u_\eps)_\eps \simh (v_\eps)_\eps$, if the following conditions
are satisfied:
\begin{itemize}
\item[(i)] For each $K\comp X$,
$\sup_{p\in X}d_h(\underline{u_\eps},\underline{v_\eps}) \to 0$ 
for some (hence every) Riemannian metric $h$ on $Y$.
\item[(ii)] $\forall k\in \N_0\ \forall m\in \N\ \forall (V,\vphi)$
chart in $X$, $\forall (W,\Psi)$ vector bundle
chart in $F$, $\forall L\comp V\ \forall L'\comp W
\ \exists \eps_1>0\ \exists C>0$ such that
$$
\|D^{(k)}(\Psi\circ u_\eps\circ\vphi^{-1} - \Psi\circ v_\eps\circ\vphi^{-1})
(\vphi(p))\| \le C\eps^m
$$
for each $\eps<\eps_1$ and each $p\in
L\cap\underline{u_\eps}^{-1}(L')\cap
\underline{v_\eps}^{-1}(L')$.
\end{itemize}
\end{definition} 

Since $\pro_1\circ \Psi\circ u_\eps\circ \vphi^{-1} = \psi\circ 
\underline{u_\eps}\circ\vphi^{-1}$ (with $\pro_1: 
\psi(W)\times \R^{n'} \to \psi(W)$)
it follows that 
\ref{hybridequiv}\, (i) precisely means that 
$(\underline{u_\eps})_\eps \sim (\underline{v_\eps})_\eps$ in $\esm[X,Y]$.
Moreover, by the same methods as employed 
in \cite{gfvm},  Remarks 2.4 and 2.6 
we conclude
that both moderateness and equivalence can be formulated equivalently by merely
requiring 
\ref{hybridmod} resp.\ 
\ref{hybridequiv} 
for charts from any given
atlas of $X$ resp.\ vector bundle atlas of $F$.
\bd The space of hybrid Colombeau generalized functions from the manifold $X$ 
into the vector bundle $F$ is defined by 
$$
\gsh[X,F] := \esmh[X,F]\big/\simh\,.
$$
\end{definition}

\bt \label{comphybrid}\  
\begin{itemize}
\item[(i)] Let $u=\cl[(u_\eps)_\eps]\in \gs[X,Y]$, $v=\cl[(v_\eps)_\eps] \in 
\Ga_{\gs}(Y,F)$. Then $v\circ u := \cl[(v_\eps\circ u_\eps)_\eps]$ is a 
well-defined element of $\gsh[X,F]$.
\item[(ii)] Let $u=\cl[(u_\eps)_\eps]\in \Ga_{\gs}(X,E)$, $v=\cl[(v_\eps)_\eps] \in 
\mathrm{Hom}_{\gs}[E,F]$. Then $v\circ u := \cl[(v_\eps\circ u_\eps)_\eps]$ is a 
well-defined element of $\gsh[X,F]$.
\end{itemize}
\end{theorem}

\pr (i) We first have to show that $(v_\eps\circ u_\eps)_\eps \in \esmh[X,F]$.
Property (i) of 
\ref{hybridmod} is obvious. 
Let $(V,\vphi)$ be a chart in $X$, $L\comp V$, $(W,\Psi)$ a vector bundle 
chart in $F$, $L'\comp W$ and $p\in L$ such that $\underline{v_\eps\circ u_\eps}(p)
= u_\eps(p) \in L'$. Then  
$D^{(k)}(\Psi\circ v_\eps\circ u_\eps\circ \vphi^{-1})(\vphi(p))
= D^{(k)}((\Psi\circ v_\eps\circ \psi^{-1}) \circ (\psi\circ u_\eps\circ \vphi^{-1}))(\vphi(p))$
can immediately be estimated using moderateness of $(u_\eps)_\eps$ and $(v_\eps)_\eps$.

To show that $v\circ u$ is well-defined let $(u_\eps)_\eps \sim (u_\eps')_\eps$
Then $(\underline{v_\eps\circ u_\eps})_\eps = (u_\eps)_\eps
\sim (u_\eps')_\eps = (\underline{v_\eps\circ u_\eps'})_\eps$, which verifies
(i) of 
\ref{hybridequiv}. To show 
(ii)
let $L'$ such that $u_\eps(L)\cup u_\eps'(L) \subseteq L'$ for $\eps$ small,
cover $L'$ 
by charts $(W_j,\psi_j)$ ($1\le j\le r$) and write $L'=\bigcup_{j=1}^r L_j'$
with $L_j'\comp \overline{W_j'}\comp W_j$ where $W_j'$ is an open
neighborhood of $L_j'$.
Then for $\eps$ small and $p\in L$ with $u_\eps(p), u_\eps'(p) \in L'$ 
there exists $j$ with $u_\eps(p),\, u_\eps'(p)\in W_j'$. 
Hence  
\beas
\lefteqn{ D^{(k)}(\Psi\circ v_\eps\circ u_\eps\circ \vphi^{-1} - 
\Psi\circ v_\eps\circ u_\eps'\circ \vphi^{-1})(\vphi(p))}
\\&&\hphantom{mm}
   =D^{(k)}((\Psi\circ \Psi_j^{-1}) \circ (\Psi_j
\circ v_\eps \circ \psi_j^{-1}) \circ (\psi_j\circ u_\eps\circ \vphi^{-1})
\\&&\hphantom{mmmm}
   -(\Psi\circ \Psi_j^{-1}) \circ (\Psi_j
\circ v_\eps \circ \psi_j^{-1}) \circ (\psi_j\circ u_\eps'\circ \vphi^{-1}))
(\vphi(p)).
\eeas
Since the norm of any derivative of each
$\Psi_{j} \circ v_\eps \circ \psi_{j}^{-1}$ is bounded by some
inverse power of $\eps$ uniformly on $W_{j}'$ ($1\le j\le r$), $(v_\eps\circ u_\eps)
\simh (v_\eps\circ u_\eps')_\eps$ follows from the above equality and 
$(u_\eps)_\eps \sim (u_\eps')_\eps$ using \cite{gfvm}, Lemma 2.5.

Finally, let $(v_\eps - v_\eps')_\eps \in \Ga_{\ns}(Y,F)$. We have to show
that
$(v_\eps\circ u_\eps)_\eps \simh (v_\eps'\circ u_\eps)_\eps$. In this case,
(i) of 
\ref{hybridequiv} is satisfied trivially. Concerning (ii), let
$L\comp V$, $L'\comp W$, $p\in L\cap u_\eps^{-1}(L')$. Then
$\|D^{(k)}(\Psi\circ(v_\eps-v_\eps')\circ\psi^{-1})\circ(\psi\circ
u_\eps\circ
\vphi^{-1})\|$ can be estimated using the $\ns$-bounds for
$(v_\eps-v_\eps')_\eps$
on $L'$ and moderateness of $(u_\eps)_\eps$, yielding the claim.

(ii) In this part of the proof we only record the general structure of the terms
to be estimated and do not embark on the topological arguments (which, anyways,
are of a simpler nature than in (i)). Since in this case
$\underline{v_\eps\circ u_\eps} = \underline{v_\eps}$, (i) of 
\ref{hybridmod}  is again obvious. Moderateness now follows by estimating
terms of the form $D^{(k)}((\Psi\circ v_\eps\circ \Phi^{-1})\circ (\Phi\circ u_\eps
\circ \vphi^{-1}))$. To see that the composition is well-defined, suppose first
that $(u_\eps)_\eps - (u_\eps')_\eps \in \Ga_{\ns}(X,E)$. Now we note that
(writing $\Phi\circ u_\eps \circ \vphi^{-1} = \big(x \mapsto
(x,u_{\eps \Phi}^{(2)}(x))\big)$)
\begin{equation} \label{desdoda}
(\Psi\circ v_\eps \circ \Phi^{-1})\circ (\Phi\circ u_\eps\circ \vphi^{-1})(x)
= (v_{\eps \Psi\Phi}^{(1)}(x), v_{\eps \Psi\Phi}^{(2)}(x)u_{\eps
\Phi}^{(2)}(x)).
\end{equation}
Hence we have to estimate $v_{\eps \Psi\Phi}^{(2)}(x)(u_{\eps \Phi}^{(2)}(x)
- {u'}_{\eps \Phi}^{(2)}(x))$ which is immediate from our assumption on $(u_\eps)_\eps$,
$(u_\eps')_\eps$. Finally, let $(v_\eps)_\eps\sim_{vb} (v_\eps')_\eps$. Then to show
equivalence of $(v_\eps\circ u_\eps)_\eps$ and $(v_\eps'\circ u_\eps)_\eps$,
by (\ref{desdoda}) it suffices to consider $(v_{\eps \Psi\Phi}^{(2)}(x) -
{v'}_{\eps \Psi\Phi}^{(2)}(x)) u_{\eps \Phi}^{(2)}(x)$, so the claim follows. 
\ep

\vskip6pt
Using $\gsh$ we now introduce the notion of generalized sections along generalized
maps, and, in particular, of generalized vector fields on
generalized maps.
\bd \label{vfalonggenmap} 
For $u \in \gs[X,Y]$ we denote by $(\gsh[X,F])(u)$ 
the set
$\{\xi\in \gsh[X,F]\mid \underline{\xi} = u\}$. In particular,
a generalized vector field $\xi$ on $u \in \gs[X,Y]$ is an 
element of
$$
\mathfrak{X}_{\gs}(u) := \{\xi\in\gsh[X,TY] \mid 
\underline{\xi} = u \}.
$$
\end{definition}

\section{Generalized pseudo-Riemannian Geometry}\label{genprg}

The aim of this section is to initiate a study of pseudo-Riemannian
geometry in the present setting. We start by introducing the notion of
a generalized connection and its Christoffel symbols.
\bd\label{gcon}\ 
\begin{itemize} 
\item[(i)] A {\em generalized connection $\hat D$} on a manifold $X$ is a map
${\gs}^1_0(X)\times{\gs}^1_0(X)\to{\gs}^1_0(X)$ satisfying 
\begin{itemize}
\item[(D1)] $\hat D_\xi \eta$ is $\gR$-linear in $\eta$.
\item[(D2)] $\hat D_\xi \eta$ is $\gs(X)$-linear in $\xi$.
\item[(D3)] $\hat D_\xi(u\eta)=u\,\hat D_\xi \eta+\xi(u)\eta$ for all $u\in\gs(X)$.
\end{itemize}
\item[(ii)] Let $(V_\al,\psi_\al)$ be a chart on $X$ with coordinates $x^i$.
We define the {\em generalized Christoffel symbols} for this chart to be the
$n^3$ functions $\hat \Ga^k_{ij}\in\gs(V_\al)$ given by
\[ 
        \hat D_{\pa_i}\pa_j\,=\,\sum\limits_k\hat \Ga^k_{ij}\,\pa_k\quad1\leq i,j\leq n\,.
\]
\end{itemize}
\end{definition}

Since ${\mathcal C}^\infty(X)$ is a submodule of $\gs(X)$ and the sheaf $\gs(X)$ is fine,
(D2) and (D3) in particular imply localizability of any generalized connection
with respect to its arguments.

We are now in the position to prove the ``Fundamental Lemma of (pseudo)-Riemannian 
Geometry'' in the present setting.
\bt
Let $(X,\hat g)$ be a generalized pseudo-Riemannian manifold.\,Then there exists a unique
generalized connection $\hat D$ such that
\begin{itemize}
\item[(D4)] $[\xi,\eta]=\hat D_\xi \eta-\hat D_\eta\xi$ and
\item[(D5)] $\xi\,\hat g(\eta,\zeta)=\hat g(\hat D_\xi \eta,\zeta)+ \hat g(\eta,\hat D_\xi \zeta)$
\end{itemize}
hold for all $\xi,\,\eta,\,\zeta$ in ${\gs}^1_0(X)$. 
$\hat D$ is called {\em generalized Levi-Civita connection} of $X$
and characterized by the Koszul formula
\beq \label{koszul}
\begin{array}{rcl}
        2\hat g(\hat D_\xi \eta,\zeta)
        &=& \xi \hat g (\eta,\zeta)+ \eta \hat g(\zeta,\xi) - \zeta \hat g(\xi,\eta)\\[.5em]
        && -\hat g(\xi,[\eta,\zeta])+\hat g(\eta,[\zeta,\xi])+ \hat g(\zeta,[\xi,\eta])\,.
\end{array}
\eeq
\end{theorem}

\pr Assume $\hat D$ to be a generalized connection additionally satisfying (D4) and (D5).
As in the classical proof (see e.g., \cite{oneill}, 
\S3, theorem 11) using the latter two
properties one shows that equation~(\ref{koszul}) is satisfied and by the 
injectivity of the map in \ref{mlemma} (ii), uniqueness follows.
\newline
To show existence define $F(\xi,\eta,\zeta)$ to be one half the right hand side of~(\ref{koszul}). Then
for fixed $\eta,\zeta$ the function $\xi\mapsto F(\xi,\eta,\zeta)$ is $\gs(X)$-linear, hence defines
a generalized one-form (using (\ref{gsrmodule})). 
Again by \ref{mlemma} (ii) there exists a unique generalized vector field 
metrically equivalent to this one-form which we may call $\hat D_\xi \eta$.
Now it is easy to derive (D1)-(D5) along the lines of the classical proof just using
the bilinearity of $g$ and the standard properties of the Lie bracket 
(cf.\ \cite{ndg}, the remark following Def.\ 10).
\ep

As in the classical case from the torsion-free condition (i.e., (D4)) 
we immediately infer
the symmetry of the Christoffel symbols of the Levi-Civita connection in the lower 
pair of indices. Moreover, from (D3) and the Koszul formula~(\ref{koszul}) we 
derive (analogously to the classical case) the following
\bp Given a chart as in \ref{gcon} (ii) we have for the 
generalized Levi-Civita connection $\hat D$ of $(X,\hat g)$ and any vector 
field $\xi\in{\gs}^1_0(X)$
\[ 
        \hat D_{\pa_i}(\xi^j\pa_j)
        \,=\,\left(\frac{\pa\,\xi^k}{\pa x^i}
        +\hat \Ga^k_{ij}\,\xi^j\right)\,\pa_k\,.
\]
Moreover, the generalized Christoffel symbols are given by
\[
        \hat \Ga^k_{ij}\,=\,\frac{1}{2}\,\hat g^{km}\,
        \left(\frac{\pa \hat g_{jm}}{\pa x^i}+\frac{\pa \hat g_{im}}{\pa x^j}
        -\frac{\pa \hat g_{ij}}{\pa x^m}\right)\,.
\]
\end{proposition}

In particular we see that we could equivalently have introduced the generalized 
Christoffel symbols of a generalized metric by demanding the classical formula 
on the level of representatives. 

To state the next result concerning consistency 
properties of generalized connections resp.\ generalized Christoffel
symbols with respect to their classical counterparts we need to define the action of 
a classical (smooth) connection $D$ on generalized vector fields $\xi$, $\eta$. 
This is done by setting
\[
        D_\xi \eta\,:=\,\cl[(D_{\xi_\eps}\eta_\eps)_\eps]\,,
\]
which is easily seen to be independent of the representatives chosen for
$\xi$ and $\eta$.
From the local formulae in the above proposition and \cite{ndg} Prop.\ 3 we conclude
\bp \label{5.4} 
Let $(X,\hat g)$ be a generalized pseudo-Riemannian manifold.
\begin{itemize}
\item[(i)] If $\hat g_{ab}=\Sigma(g_{ab})$ where $g_{ab}$ is a classical smooth pseudo-Riemannian metric
then we have, in any chart, $\hat \Ga^i_{jk}=\sigma(\Ga^i_{jk})$ (with $\Ga^i_{jk}$ denoting the 
Christoffel Symbols of $g$). Hence for all $\xi, \eta\in{\gs}^1_0(X)$
\[ \hat D_\xi \eta\,=\,D_\xi \eta\,.\]
\item[(ii)] If $\hat g_{ab}\approx_\infty g_{ab}$, $g_{ab}$ a classical smooth metric,
$\xi, \eta\in{\gs}^1_0(X)$ and $\xi\approx_\infty\zeta\in\T^1_0(X)$, 
$\eta\approx\nu\in\D'^1_0(X)\ ($or vice versa, i.e.,  $\xi\approx\zeta\in\D'^1_0(X)$, 
$\eta\approx_\infty\nu\in\T^1_0(X)\,)$ then
\[ \hat D_\xi \eta\,\approx D_\zeta\nu\,.\]
\item[(iii)] Let $\hat g_{ab}\approx_k g_{ab}$, $g_{ab}$ a classical $\CC^k$-metric,  
then, in any chart, $\hat\Ga^i_{jk}\,\approx_{k-1}\,\Ga^i_{jk}$.
If in addition $\xi$, $\eta\in{\gs}^1_0(X)$, $\xi\approx_{k-1}\zeta\in\Ga^{k-1}(X,TX)$ 
and $\eta\approx_k\nu\in\Ga^k(X,TX)$ then 
\[ \hat D_\xi \eta\,\approx_{k-1}D_\zeta\nu\,.\]
\end{itemize}
\end{proposition}

Our next aim is to define the induced covariant derivative of a generalized metric on 
a generalized curve. 
Let $J\subseteq \R$ be an interval and $\gamma \in \gs[J,X]$. Let $\hat g\in {\gs}^0_2(X)$ be
a generalized metric. For any $K\comp J$ there exists $\eps_0>0$ and $K'\comp X$ such
that $\gamma_\eps(K)\subseteq K'$ for $\eps<\eps_0$. According to \ref{mainmetric} we
may choose a representative $(\hat g_\eps)_\eps$ of $\hat g$ such that each $\hat g_\eps$ is
a pseudo-Riemannian metric in a neighborhood of $K'$. 
Let $\xi = \cl[(\xi_\eps)_\eps]\in \mathfrak{X}_{\gs}(\gamma)$. For each fixed small $\eps$ 
we let $\xi_\eps'$ be the induced covariant derivative of $\xi_\eps$ on $\gamma_\eps$ 
with respect to $\hat g_\eps$. 
\bd
We call $\xi':= \cl[(\xi_\eps')_\eps]\in \mathfrak{X}_{\gs}(\gamma)$ the {\em induced
covariant derivative} of $\xi$ on $\gamma$ with respect to $\hat g$.
\end{definition}

For this definition to make sense
we have to show that $\xi'$ is independent of the chosen representatives 
$(\gamma_\eps)_\eps$, $(\xi_\eps)_\eps$
and $(\hat g_\eps)_\eps$. To this end we note that for fixed $\eps$ the local form of
$\xi_\eps'$ is given by (with $\hat\Ga^k_{\eps ij}$ the Christoffel symbols of $\hat g_\eps$)
\begin{equation} \label{indcovderloc}
\xi_\eps' = \sum_k \left(\frac{d\xi^k_\eps}{dt}  + 
\sum_{i,j} \hat\Ga^k_{\eps ij} \frac{d\gamma_\eps^i}{dt}\xi^j_\eps \right)\pa_k\,.
\end{equation}
From this on the one hand we conclude that $(\xi_\eps')_\eps$ is indeed moderate 
and on the other hand it is straightforward to check that
choosing different representatives for $\xi$, $\gamma$ or
$\hat g$ does not change the class of $\xi'$ in $\mathfrak{X}_{\gs}(\gamma)$. Finally, we 
conclude that the above restriction of $\gamma$ to relatively compact subintervals of
$J$ can be overcome by ``patching together'' the representatives of $\xi'$ obtained
for a covering of $J$ by relatively compact subintervals. In fact, these partially
defined generalized functions coincide on overlapping intervals again due
to the explicit local form (\ref{indcovderloc}).

The main properties of the induced covariant derivative are collected in the following
result.
\bp
Let $\hat g$ be a generalized metric on $X$ with Levi Civita connection $\hat D$ and
let $\gamma \in \gs[J,X]$. Then 
\begin{itemize}
\item[(i)] $(\tilde r \xi_1 + \tilde s \xi_2)' = \tilde r \xi_1' + \tilde s \xi_2'$ 
$\qquad (\tilde r,\, \tilde s \in \ks, \, \xi_1,\, \xi_2 \in \mathfrak{X}_{\gs}(\gamma))$.
\item[(ii)] $(u\xi)' = \frac{du}{dt}\xi + u\xi'$ $\qquad (u\in \gs(J),\,
\xi \in \mathfrak{X}_{\gs}(\gamma))$.
\item[(iii)] $(\xi\circ \gamma)' = \hat D_{\gamma'(.)}\xi$ 
\hspace{1em} in $\mathfrak{X}_{\gs}(\gamma)$ $(\xi \in {\gs}^1_0(X))$.
\end{itemize}
\end{proposition}

\pr The composition $\xi\circ \gamma$ in (iii) 
is well-defined by \ref{comphybrid} (i).
Moreover, the right hand side of (iii) exists since $\xi\to \hat D_\xi\eta$ is a (vector valued)
generalized tensor field (cf.\ (D2) in \ref{gcon}).
All the claimed identities now follow directly from the local form 
(\ref{indcovderloc}) of the induced covariant derivative. 
\ep

On $J$ we consider the section $s:= t \mapsto \left.\frac{d}{dt}\right|_{t} \equiv (t,1)$.
Applying \cite{gfvm}, Def.\ 3.3
to a generalized curve $\ga\in\gs[J,X]$, it follows from \ref{comphybrid} (ii) 
that its {\em velocity
vector field} defined by $\ga':= T\ga \circ s$ 
is a well-defined element of $\mathfrak{X}_{\gs}(\gamma)$. Then $\ga''$
is defined as the induced covariant derivative of $\ga'$ on $\ga$.
\bd A geodesic in a generalized pseudo-Riemannian manifold is a curve
$\ga\in\gs[J,X]$ ($J\subseteq\R$) satisfying $\ga''=0$.
\end{definition}

By (\ref{indcovderloc}) $\ga$ is a geodesic in $(X,\hat g)$ if and only if the usual
local formula holds, i.e., if and only if 
\begin{equation}\label{geo}
  \Big[\,\Big(\frac{d^2\gamma_\eps^k}{dt^2}
  +\sum_{i,j}\hat\Ga^k_{\eps ij}\frac{\gamma_\eps^i}{dt}\frac{\gamma_\eps^j}{dt}\Big)_\eps\,\Big]
  =0\quad\mbox{ in }\mathfrak{X}_\gs(\ga)\,.
\end{equation}

\section{Applications to general relativity}\label{gengr}
Before introducing the generalized curvature tensor and its contractions
entering Einstein's equations we briefly comment on the shortcomings of 
classical methods and in particular of the linear distributional geometry 
(as introduced e.g., in \cite{marsden,parker}) especially in the
context of general relativity. 
While successfully used within linear field theories (e.g., point charges in 
electrodynamics), applications of distributional methods to general relativistic 
problems have been rare in the literature;
the source of all 
difficulties of course is the nonlinearity of the field equations. 

On the one hand, in the context of (local) existence and uniqueness theorems 
for Einstein's equations one is bound to work with Sobolev spaces of high enough
regularity so that the metric belongs to some classical function algebra.
Roughly speaking, the classical local existence theorems guarantee the existence 
of a unique solution $g$ to the equations (formulated as an initial value
problem along a spacelike hypersurface $\Sigma$ for a Riemannian metric $h$ and its extrinsic 
curvature $K$ satisfying the constraint equations) with the metric
$g\in\CC^0([0,T);H^s(\Sigma)))\cap \CC^1([0,T);H^{s-1}(\Sigma))$ for 
data $h,K$ chosen in $H^{s-1}(\Sigma)$ provided that  $s> 4$ (cf. \cite{klainerman1,rendall}). 
In fact the optimal local existence result in the case of asymptotic flatness actually only 
requires $s > 5/2$ and recent developments by Klainerman und Rodnianski \cite{klainerman2} aim at 
further improving this to the bound $s > 2$. 

On the other hand, when dealing with special (constructive) solutions 
of Einstein's equations one has to cope with the problem of computing the
curvature of a given metric of low differentiability (in particular, distributional)
which again is problematic due to the nonlinearities involved. Within the $\D'$-framework
it is nevertheless possible to consistently describe  sources of the gravitational field concentrated 
(i.e., the energy-momentum tensor supported) on a submanifold of codimension one 
in spacetime (so-called thin shells, cf.\ \cite{israel1}). 
However, in a classical paper \cite{gt} Geroch and Traschen have shown that 
within classical (linear) distribution theory gravitating sources confined to 
a submanifold of codimension greater than one in 
spacetime (hence, in particular, such interesting objects as cosmic strings) are 
excluded from a mathematically rigorous and at the same time physically sensible 
description. By the latter we mean the existence of an appropriate notion of 
convergence of metrics which ensures the convergence of the respective curvature 
tensors. 

Here we are going to introduce a setting that is primarily intended to cope with the
latter situation described above which at the same time is mathematically rigorous and
physically sensible.

We start by defining the generalized Riemann, Ricci, scalar and 
Einstein curvature from an 
invariant point of view. It is then clear that all the classical formulae will hold on 
the level of representatives, i.e., all the symmetry properties of the respective 
classical tensor fields carry over to the new setting. Moreover, the Bianchi identities 
hold in the generalized sense.
\bd \label{cquant}
Let $(X,\hat g)$ be a generalized pseudo-Riemannian manifold with Levi-Civita 
connection $\hat D$.
\begin{itemize}
\item[(i)] The {\em generalized Riemannian curvature tensor}
$ \hat R_{abc}\,^d\in{\gs}^1_3(X)$ is defined by
\[       \hat R_{\xi,\eta}\zeta\,:=\, \hat D_{[\xi,\eta]} \zeta-[\hat D_\xi,\hat D_\eta]\zeta\,.
\]
\item[(ii)] 
We define the {\em generalized Ricci curvature tensor} $ \hat R_{ab}\in{\gs}^0_2(X)$ 
by the usual contraction of the generalized Riemann tensor
\[
         \hat R_{ab}\,:=\, \hat R_{cab}\,^c\,.
\]
\item[(iii)] 
The {\em generalized curvature scalar (or Ricci scalar)} $ \hat R\in\gs(X)$ is defined
by the usual contraction of the generalized Ricci tensor
\[
         \hat R\,:=\, \hat g^{ab} \hat R_{ab}\,. 
\]
\item[(iv)] 
Finally we define the {\em generalized Einstein tensor} $\hat G_{ab}
\in{\gs}^0_2(X)$ by
\[
        \hat G_{ab}\,:=\,\hat R_{ab}-\frac{1}{2}\,\hat R\,\hat g_{ab}\,.
\]
\end{itemize}
\end{definition}

The framework developed above opens the gate to a wide range of applications
in general relativity. 
\ref{gmdef}  is capable of modelling a large class
of singular spacetimes while at the same time its (generalized) curvature quantities 
simply may be calculated by the usual coordinate formulae.
Hence we are in a position to mathematically rigorously formulate 
Einstein's equations for generalized metrics. Moreover we have at our disposal 
several theorems (which essentially are rooted in \cite{ndg}, Prop. 3) guaranteeing 
consistency with respect to linear distributional geometry resp.\ the smooth setting. 

\bt \label{6.2}
Let $(X,\hat g)$ a generalized pseudo-Riemannian manifold with 
$\hat g_{ab}\approx_k g_{ab}$.
Then all the generalized curvature quantities defined above are 
$\CC^{k-2}$-associa\-ted with their classical counterparts. 
\end{theorem}

In particular, if a generalized metric $\hat g_{ab}$ is $\CC^2$-associated 
to a  vacuum solution of Einstein's equations then we have for the generalized Ricci tensor
\[ 
    \hat R_{ab}\,\approx_0 R_{ab}\,=0\,.
\]
Hence $\hat R_{ab}$ satisfies the vacuum Einstein equations in the sense of 
$0$-association (cf.\ the remarks after Prop.\  18 in the revised version 
of \cite{sotonTF}).
\ms
Generally speaking, whenever we encounter a spacetime metric of low differentiability 
in general relativity we may proceed 
along the lines of the following blueprint to obtain a mathematically and physically 
satisfactory description of the singular spacetime geometry: 
first we have to transfer the classically singular metric to 
a generalized one. This may be done by some ``canonical'' smoothing or by some 
other physically motivated regularization (see also the remarks on 
nonlinear modelling preceding Def.\ 2 in \cite{ndg}). 
Once the generalized setting has been entered,
all curvature quantities may be calculated simply using componentwise classical
calculus. All classical concepts literally carry over to the new framework and one
may treat e.g., the Ricci tensor, geodesics, geodesic deviation, etc.\ within this 
nonlinear distributional geometry. Finally one may use the concept of association
to return to the distributional or $\CC^k$-level for the purpose of interpretation.

This program has been carried out for a conical metric (representing a cosmic string) 
by Clarke, Vickers and Wilson (see \cite{clarke,invc,wilson} for a treatment 
in the full 
setting of Colombeau's construction) rigorously assigning to it a distributional 
curvature and (via the field equations) the heuristically expected energy-momentum 
tensor. Moreover, the nonlinear generalized function setting was used
in \cite{hb5,Balkerr} to calculate the energy momentum tensor
of the extended Kerr geometry as well as in \cite{mark} to unify several distributional 
approaches to the Schwarzschild geometry. Finally, 
a complete distributional description
of impulsive pp-wave spacetimes was achieved  
in \cite{geo,geo2,penrose}.
\ms
An in-depth analysis of the interrelations between the setting introduced here and 
the one of Geroch and Traschen (\cite{gt}) as well as a study of spherical impulsive
gravitational waves in this framework is 
the subject of ongoing research.

We close this work by discussing the geodesic equation of
impulsive pp-waves in the present setting. In \cite{geo,geo2,penrose} the geometry
of an impulsive pp-wave was described by the following generalized
line element
(cf.\ \ref{mex} (iii))
\beq\label{pp}
\hat{ds}^2=f(x,y)D(u)du^2-dudv+dx^2+dy^2\,.
\eeq
Here the spacetime manifold $X$ is taken to be $\R^4$ and $D$ is a generalized 
delta function, i.e., $D=\cl[(\de_\eps)_\eps]$, with $\int\de_\eps\to 1$ the support 
of $\de_\eps$ shrinking to $0$ and $\de_\eps$ uniformly bounded in $L^1$ for small 
$\eps$. (cf. \cite{MObook}, Def.\ 7.1). 
Physically this spacetime describes a 
gravitational impulse located at the null-hypersurface $u=0$ in Minkowski space; 
the curvature vanishes everywhere but on the impulse.

(\ref{pp}) provides an example of a general regularization procedure suggested by the 
physical situation to be modelled (cf.\ the discussion in Section \ref{intro}). In fact, 
generalized delta functions in the above sense provide a very general class of 
regularizations of delta-type singularities. (\ref{pp}) complies with 
viewing the singular metric itself as an impulsive limit of sandwich waves of 
infinitely short duration.

We have to solve the system (\ref{geo}) for the metric (\ref{pp}).
Due to the special form of the metric the first equation (i.e., $k=0$) is trivial, so that
(using coordinates $u,v$ and $x^i=(x,y)$ as in (\ref{pp})) $u$ may be used as an 
affine parameter along the geodesics and system (\ref{geo}) reduces to 
\beas
        \ddot v_\eps(u)&=&f(x_\eps^j(u))\,\dot\de_\eps(u)
              \,+\,2\,\pa_if(x_\eps^j(u))\,\dot x_\eps^i\,(u)\de_\eps(u)\nn \\
        \ddot x_\eps^i\,(u)&=&\frac{1}{2}\,\pa_i\,f(x_\eps^j(u))\,\,\de_\eps(u)
\eeas
This system of nonlinear ODEs was shown to be uniquely solvable in
$\gs(\R)^3$ (cf.\ \cite{geo2}, Th.\ 1) for given initial conditions 
$(v_\eps,x_\eps^i)(-1)=(v_0,x^i_0)$ and $(\dot v_\eps,\dot x^i_\eps)(-1)=
(\dot v_0,\dot x^i_0)$. The resulting geodesic is an element of $\gs[\R,X]$. 
Moreover, the generalized solution (by \cite{geo2}, Th.\ 3)
possesses the following distributional shadow 
\beas
        x^i &\approx & x_0^i + \dot x_0^i (1+u) + \frac{1}{2}\pa_i f(x_0^j +
                \dot x_0^j) u_+ \nn\\
        v\, &\approx & v_0 + \dot v_0 (1+u) + f(x_0^j +
                \dot x_0^j) H(u) \\ 
                && + \pa_if(x_0^j + \dot x_0^i) \left(\dot x_0^i +
        \frac{1}{4}\pa_i f(x_0^j +  \dot x_0^j)\right) u_+\,\,,
\eeas
where $H$ denotes the Heaviside function and $u_+=H(u)u$.
Hence the macroscopic aspect of the generalized geodesics displays the physically
sensible behavior of the geodesics being refracted broken straight lines. 
It should be noted, though, that the explicit calculation of distributional limits
in fact is based on the fact that the underlying manifold is given as $\R^4$ since 
in linear distribution theory (contrary to the Colombeau setting) there is no concept
of generalized functions valued in a manifold.

Finally we would like to emphasize that the framework developed in 
the previous sections for the first time allows 
a comprehensive and consistent interpretation of the 
calculations given, e.g.,  in \cite{clarke, mark, geo2, penrose}.

{\bf Acknowledgments.} We would like to thank Michael Grosser and Michael
Oberguggenberger for numerous discussions that importantly contributed to the
final form of the manuscript. Also, we are indebted to the referees for 
several helpful comments.


\begin{thebibliography}{{Kun}99b}

\bibitem{AJ}
{Aragona, J., Juriaans, S.\ O.}
\newblock Some structural properties of the topological ring of {C}olombeau's
  generalized numbers.
\newblock {\em Comm.\ Algebra}, {\bf 29} (5), 2201--2230, 2001.

\bibitem{hb5}
{Balasin, H.}
\newblock Distributional energy-momentum tensor of the extended {K}err
  geometry.
\newblock {\em Class.\ Quant.\ Grav.}, pages 3353--3362, 1997.

\bibitem{Balkerr}
{Balasin, H.}
\newblock Distributional aspects of general relativity: The example of the
  energy-momentum tensor of the extended {K}err-geometry.
\newblock In {Grosser, M., Hörmann, G., Kunzinger, M., Oberguggenberger, M.},
  editor, {\em Nonlinear Theory of Generalized Functions}, volume 401 of {\em
  CRC Research Notes}, pages 231--239, Boca Raton, 1999. CRC Press.

\bibitem{bhatia}
{Bhatia, R.}
\newblock {\em Perturbation Bounds for Matrix Eigenvalues}, volume~{\bf 162} of
  {\em Pitman Research Notes in Mathematics}.
\newblock Longman, Harlow, U.K., 1987.

\bibitem{bourbaki-algebra}
{Bourbaki, N.}
\newblock {\em Algebra I, Chapters 1--3}.
\newblock Elements of Mathematics. Hermann, Addison-Wesley, Paris,
  Massachusetts, 1974.

\bibitem{clarke}
{Clarke, C.~J.~S., Vickers, J.~A., Wilson, J.~P.}
\newblock Generalised functions and distributional curvature of cosmic strings.
\newblock {\em Class.~Quant.~Grav.}, {\bf 13}, 1996.

\bibitem{c1}
{Colombeau, J.~F.}
\newblock {\em New Generalized Functions and Multiplication of Distributions}.
\newblock North Holland, Amsterdam, 1984.

\bibitem{c2}
{Colombeau, J.~F.}
\newblock {\em Elementary Introduction to New Generalized Functions}.
\newblock North Holland, Amsterdam, 1985.

\bibitem{DKP}
{Dapi\'c, N., Kunzinger, M., Pilipovi\'c, S.}
\newblock Symmetry group analysis of weak solutions.
\newblock {\em Proc.\ London Math.\ Soc., to appear}, 2002.

\bibitem{RD}
{De Roever, J.~W., Damsma, M.}
\newblock Colombeau algebras on a {${\mathcal C}^\infty$}-manifold.
\newblock {\em {Indag.~Mathem., N.S.}}, {\bf 2}(3), 1991.

\bibitem{gt}
{Geroch, R., Traschen, J.}
\newblock Strings and other distributional sources in general relativity.
\newblock {\em Phys.~Rev.~D}, {\bf 36}(4):1017--1031, 1987.

\bibitem{found}
{Grosser, M., Farkas, E., Kunzinger, M., Steinbauer, R.}
\newblock On the foundations of nonlinear generalized functions {I}, {II}.
\newblock {\em Mem.~Am.~Math.~Soc.} {\bf 153}, 2001.

\bibitem{vim}
{Grosser, M., Kunzinger, M., Steinbauer, R., Vickers, J.}
\newblock A global theory of algebras of generalized functions.
\newblock {\em Adv.\ Math.}, to appear, 2001.

\bibitem{mark}
{Heinzle, J.\ M., Steinbauer, R.}
\newblock Remarks on the distributional Schwarzschild geometry.
\newblock {\em J.~Math.~Phys., to appear}, 2001.

\bibitem{israel1}
{Israel, W.}
\newblock Singular hypersurfaces and thin shells in general relativity.
\newblock {\em Nouv.~Cim.}, {\bf 44B}(1):1--14, 1966.

\bibitem{kamleh}
{Kamleh, W.}
\newblock Signature changing space-times and the new generalised functions.
\newblock {\em Preprint, gr-qc/0004057}, 2000.

\bibitem{klainerman1}
{Klainerman, S., Nicol\'o, F.}
\newblock On local and global aspects of the Cauchy problem in general relativity.
\newblock {\em Calss. Quant. Grav.}, {\bf 16}:R73--R157, 1999.

\bibitem{klainerman2}
{Klainerman, S., Rodnianski, I.}
\newblock Rough solutions of the Einstein-vacuum equations.
\newblock {\em Preprint, math.AP/0109173}, 2001. 

\bibitem{KM}
{Kriegl, A., Michor, P.~W.}
\newblock {\em The Convenient Setting of Global Analysis}, volume~{\bf 53} of
  {\em Math.~Surveys Monogr.}
\newblock Amer.~Math.~Soc., Providence, RI, 1997.

\bibitem{gfvm}
{Kunzinger, M.}
Generalized functions valued in a smooth manifold.
Submitted (available electronically at 
http://arxiv.org/abs/math.FA/ 0107051 ), 2001.

\bibitem{symm}
{Kunzinger, M., Oberguggenberger, M.}
\newblock Group analysis of differential equations and generalized functions.
\newblock {\em SIAM J.~Math.~Anal.}, {\bf 31}(6):1192--1213, 2000.


\bibitem{geo2}
{Kunzinger, M., Steinbauer, R.}
\newblock A rigorous solution concept for geodesic and geodesic deviation
  equations in impulsive gravitational waves.
\newblock {\em J.~Math.~Phys.}, {\bf 40}:1479--1489, 1999.

\bibitem{penrose}
{Kunzinger, M., Steinbauer, R.}
\newblock A note on the {P}enrose junction conditions.
\newblock {\em Class.~Quant.~Grav.}, {\bf 16}:1255--1264, 1999.


\bibitem{ndg}
{Kunzinger, M., Steinbauer, R.}
Foundations of a nonlinear distributional geometry.
{\em Acta Appl. Math.}, to appear (available electronically at 
http://arxiv.org/abs/math.FA/ 0102019), 2001.

\bibitem{Mal} {Mallios, A.} {\em Geometry of vector sheaves}. Vol. I, II. 
Mathematics and its Applications, 
{\bf 439}. Kluwer Academic Publishers, Dordrecht, 1998. 

\bibitem{MallRos0}
{Mallios, A., Rosinger, E.\ E.}
 Abstract differential geometry, differential algebras of generalized
  functions, and de Rham cohomology.
 {\em Acta Appl.\ Math.}, {\bf 55}(3):231--250, 1999.

\bibitem{MallRos}
{Mallios, A., Rosinger, E.\ E.}
 Space-time foam dense singularities and de {R}ham cohomology.
 {\em Acta Appl.\ Math.}, {\bf 67}: 59--89, 2001.

\bibitem{LT}
{Lig\c{e}za, J., Tvrdy, M.}
\newblock On systems of linear algebraic equations in the {C}olombeau algebra.
\newblock {\em Math.\ Bohem.}, {\bf 124}(1):1--14, 1999.

\bibitem{MaN}
{Mansouri, R., Nozari, K.}
\newblock A new distributional approach to signature change.
\newblock {\em Gen. Relativity Gravitation}, {\bf
  32}(2):235--269, 2000.

\bibitem{marsden}
{Marsden, J.~E.}
Generalized {H}amiltonian mechanics.
{\em Arch.\ Rat.\ Mech.\ Anal.}, {\bf 28}(4):323--361, 1968.

\bibitem{MObook}
{Oberguggenberger, M.}
\newblock {\em Multiplication of Distributions and Applications to Partial
  Differential Equations}, volume~{\bf 259} of {\em Pitman Research Notes in
  Mathematics}.
\newblock Longman, Harlow, 1992.

\bibitem{point}
{Oberguggenberger, M., Kunzinger, M.}
\newblock Characterization of {C}o\-lom\-beau generalized functions by their
  pointvalues.
\newblock {\em Math.~Nachr.}, {\bf 203}:147--157, 1999.

\bibitem{oneill}
{O'Neill, B.}
\newblock {\em Semi-Riemannian Geometry (With Applications to Relativity)}.
\newblock Academic Press, New York, 1983.

\bibitem{parker}
{Parker, P.}
\newblock Distributional geometry.
\newblock {\em J.~Math.~Phys.}, {\bf 20}(7):1423--1426, 1979.

\bibitem{penrose_rindler}
{Penrose, R., Rindler, W.}
\newblock {\em Spinors and space-time {I}}.
\newblock Cambridge University Press, 1984.

\bibitem{rendall}
{Rendall, A.}
\newblock Local and Global Existence Theorems for the Einstein Equations.
\newblock {\em Living Rev. Relativity}, {\bf 3}(1). [Online Article]:
cited on 2002-01-20, http://www.livingreviews.org/Articles/Volume3/2000-1rendall/.

\bibitem{simanca}
{Simanca, S.\ R.}
\newblock {\em Pseudo-differential Operators}, volume~{\bf 236} of {\em Pitman
Research in Notes in Mathematics}.
\newblock Longman, Harlow, 1990.

\bibitem{ultra}
{Steinbauer, R.}
\newblock The ultrarelativistic {R}eissner-{N}ordstr{\o}m field in the
  {C}olombeau algebra.
\newblock {\em J.~Math.~Phys.}, {\bf 38}:1614--1622, 1997.

\bibitem{geo}
{Steinbauer, R.}
\newblock Geodesics and geodesic deviation for impulsive gravitational waves.
\newblock {\em J.~Math.~Phys.}, {\bf 39}:2201--2212, 1998.

\bibitem{vickersESI}
{Vickers, J.~A.}
\newblock Nonlinear generalized functions in general relativity.
\newblock In {Grosser, M., Hörmann, G., Kunzinger, M., Oberguggenberger, M.},
  editor, {\em Nonlinear Theory of Generalized Functions}, volume 401 of {\em
  CRC Research Notes}, pages 275--290, Boca Raton, 1999. CRC Press.

\bibitem{sotonTF}
{Vickers, J., Wilson, J.}
A nonlinear theory of tensor distributions.
ESI-Preprint (available electronically at 
  http://www.esi.ac.at/ESI-Preprints.html), {\bf 566}, 1998.

\bibitem{invc}
{Vickers, J.~A., Wilson, J.~P.}
\newblock Invariance of the distributional curvature of the cone under smooth
  diffeomorphisms.
\newblock {\em Class.~Quantum.~Grav.}, {\bf 16}:579--588, 1999.

\bibitem{wilson}
{Wilson, J.~P.}
\newblock Distributional curvature of time dependent cosmic strings.
\newblock {\em Class.~Quantum Grav.}, {\bf 14}:3337--3351, 1997.

\end{thebibliography}
\end{document}